\documentclass[reqno]{amsproc}
\usepackage{amssymb}

\newtheorem{theorem}{Theorem}[section]
\newtheorem{lemma}[theorem]{Lemma}
\newtheorem*{sublemma}{Sublemma}

\theoremstyle{definition}

\theoremstyle{remark}
\newtheorem{remark}[theorem]{Remark}

\numberwithin{equation}{section}

\newcommand{\abs}[1]{\lvert#1\rvert}

\newcommand{\bbN}{{\mathbb N}}
\newcommand{\bbR}{{\mathbb R}}
\newcommand{\bbZ}{{\mathbb Z}}
\newcommand{\identrel}[1]{=_{#1}}
\newcommand{\restrict}{{\restriction}}
\newcommand{\nullset}{\varnothing}

\begin{document}

\title{How Many Turing Degrees are There?}

\author{Randall Dougherty}
\address{Department of Mathematics, Ohio State
University, Columbus, Ohio 43210}
\email{rld@math.ohio-state.edu}
\thanks{The first author was partially supported by NSF Grant DMS 9158092.}

\author{Alexander S. Kechris}
\address{Department of Mathematics, California Institute of Technology,
Pasadena, California 91125}
\email{kechris@caltech.edu}
\thanks{The second author was partially supported by NSF Grant DMS 9619880.}

\date{}

\begin{abstract}
A Borel equivalence relation on a Polish space is {\it countable} if all of
its
equivalence classes are countable.  Standard examples of countable Borel
equivalence relations (on the space of subsets of the integers) that occur
in recursion theory are:  recursive isomorphism, Turing equivalence,
arithmetic equivalence, etc.  There is a canonical hierarchy of complexity
of countable Borel equivalence relations imposed by the notion of Borel
reducibility.  We will survey results and conjectures concerning the
problem of identifying the place in this hierarchy of these equivalence
relations from recursion theory and also discuss some of their implications.
\end{abstract}

\maketitle

The obvious answer to the question of the title is:  {\it continuum many}.
There is however a different way of looking at this question, which leads
to some very interesting open problems in the interface of recursion theory
and descriptive set theory.  Our goal in this paper is to explain the
context in which this and related problems can be formulated, i.e., the
theory of Borel equivalence relations, and survey some of the progress
to date.

\section{Formulation of the problem}

We denote by $\equiv_T$ the {\bf Turing equivalence relation} on
${\mathcal P}(\bbN )=\{X:X\subseteq \bbN\}$, which we identify with
$2^\bbN$,
viewing sets
as characteristic functions.  (We use the standard set-theoretic
convention that $n = \{0,1,\dots,n-1\}$ for all natural numbers~$n$.)
Then $\equiv_T$ is a Borel (in fact
$\Sigma^0_3$) equivalence relation on $2^\bbN$.  We denote by
${\mathcal D}$ the quotient space $2^\bbN /(\equiv_T)$, i.e., the set of
{\bf Turing degrees}.

Now consider general Borel equivalence relations on $2^\bbN$ or even
arbitrary {\bf Polish} (separable completely metrizable) spaces.
We measure their complexity by studying the following partial (pre)order
of {\bf Borel reducibility}: if $E, F$ are Borel equivalence relations on
$X,Y$ respectively, then a {\bf Borel reduction} of $E$ into $F$ is a Borel
map $f:X\rightarrow Y$ such that
$$xEy\iff f(x)Ff(y).$$
If such an $f$ exists we say that $E$ is {\bf Borel reducible} to $F$
and denote this by
$$E\leq_BF.$$
Let also
$$E\sim_BF\iff E\leq_BF\ \&\ F\leq_BE$$
(this defines the concept of {\bf bi-reducibility}) and
$$E<_BF\iff E\leq_BF\ \&\ F\not\leq_BE.$$

Let us say that a function $f_*:X/E\rightarrow Y/F$ is {\bf Borel} if it
has a Borel lifting, i.e., there is a Borel function $f:X\rightarrow Y$
such that
$f_*([x]_E)=[f(x)]_F$ for all $x \in X$.  Then it
is clear that $E\leq_BF$ is equivalent to the assertion that there is
a {\bf Borel injection} from $X/E$ into $Y/F$, which we express by
saying that the {\bf Borel cardinality}, $\abs{E}_B$, of
$E$ is less than or equal to to that of $F$; in symbols,
$$\abs E_B\leq\abs F_B\iff E\leq_BF.$$
Then define
$$\abs E_B=\abs F_B\iff E\sim_BF,$$
i.e., $X/E,\ Y/F$ have the same Borel cardinality, and
$$\abs E_B<\abs F_B\iff E<_BF,$$
i.e., $X/E$ has (strictly) smaller Borel cardinality then $Y/F$.

We are now ready to formulate our problem as follows, where, by abusing
notation, we write below $\abs{\mathcal D}_B$ instead of $\abs{\equiv_T}_B$
and call this the Borel
cardinality of $\mathcal D$, instead of $\equiv_T$:

\medskip

\noindent{\bf Question}:  What is the Borel cardinality, $\abs{\mathcal
D}_B$,
of the set of Turing degrees $\mathcal D$?
\medskip

If we denote the classical (Cantor) cardinality of $\mathcal D$ by
$\abs{\mathcal D}$, then we have $\abs{\mathcal D} =\abs \bbR$.
However, it is not hard to see that the Borel cardinality
of $\mathcal D$ is bigger than that of the continuum.
Let $\identrel X$ be the identity relation on the Polish space~$X$.
So $\abs{\identrel\bbR}_B$ is the Borel cardinality which naturally
represents the classical cardinality of the continuum.
\medskip

\noindent {\bf Fact}.  $(\equiv_T)>_B(\identrel\bbR)$.

\begin{proof}
It is standard that there is a perfect set of
pairwise Turing incomparable subsets of $\bbN$, so $(\identrel\bbR)
\leq_B(\equiv_T)$.
If on the other hand $f:2^\bbN\rightarrow\bbR$ is Borel and
Turing-invariant, i.e., $x\equiv_Ty\implies f(x)=f(y)$, then for each
Borel set $A\subseteq\bbR$, $f^{-1}(A)$ is a Turing-invariant Borel
subset of $2^\bbN$, so it has measure 0 or 1.  It follows that, for each
$n$, the $n$th digit in the decimal expansion of $f(x)$ is fixed on a set of
measure 1. So there is
a Turing-invariant Borel set of measure 1 on which $f$ is constant,
therefore
$f$ cannot be a reduction of $\equiv_T$ into $\identrel\bbR$.
Thus $(\equiv_T)\not\leq_B(\identrel\bbR)$.\end{proof}

We now have our question but it is not clear yet what kind of answer we
should
expect.  In what sense can we hope to compute $\abs{\mathcal D}_B$?
To understand this, we have to dig a little deeper into the theory
of Borel equivalence relations.

For our purposes, a crucial property of the Turing equivalence relation
is that it has countable equivalence classes.  In general, we call
a Borel equivalence relation {\bf countable} if every one of its
classes is countable.  We will next review some basic facts of the theory
of countable Borel equivalence relations, for which we refer the reader to
the papers Kechris \cite{K2}, Dougherty-Jackson-Kechris \cite{DJK},
Jackson-Kechris-Louveau \cite{JKL}, Kechris \cite{K1}, and Adams-Kechris
\cite{AK}.
\medskip

(i) (Feldman-Moore \cite{FM}) Every countable Borel equivalence relation
is generated by a Borel action of a countable group.

More precisely, given a countable Borel equivalence $E$ on a Polish
space $X$, there is a countable group $G$ and a Borel action $(g,x)\mapsto
g\cdot x$ of $G$ on $X$ such that, if $E^X_G$ is defined by
$$xE^X_Gy\iff\exists g{\in}G\,\,(g\cdot x=y),$$
then $E=E^X_G$.

In particular, $\equiv_T$ is given by a Borel action of a countable group
on $2^\bbN$.  It seems like an interesting, but somewhat vague, question to
find out whether one can obtain such a representation that has some
recursion
theoretic significance.

\begin{remark}
Using the Feldman-Moore theorem and related facts, within
a Schr\"oder-Bernstein
argument, one can show that, for countable Borel equivalence
relations $E$ and~$F$,  $E\sim_BF$ is equivalent to the existence of a
Borel bijection of $X/E$ with $Y/F$.\end{remark}
\medskip

(ii) There is a {\bf universal} countable Borel equivalence relation, in
the sense of~$\leq_B$.

That is, there is a countable Borel equivalence relation $E$ such that, for
any countable Borel equivalence relation $F$, we have $F\leq_BE$.  This~$E$
is clearly unique, up to $\sim_B$, and denoted by $E_\infty$.

An example of a universal countable Borel equivalence is given by the
orbit equivalence relation of the shift action of $F_2$, the free group
on two generators, on~$2^{F_2}$ given by
$$g\cdot x(h)=x(g^{-1}h),\qquad g,h\in F_2,\quad x\in 2^{F_2}.$$

(iii) There is a smallest, in the sense of $\leq_B$, countable Borel
equivalence relation on uncountable Polish spaces, namely $\identrel\bbR$.

So for every countable Borel equivalence relation $E$ on an uncountable
Polish space, we have $(\identrel\bbR)\leq_BE$.  If
$(\identrel\bbR)\sim_BE$,
we say that $E$
is {\bf smooth}.  For example, $\equiv_T$ is not smooth.  Another example
of a non-smooth countable Borel equivalence is the following one,
defined on $2^\bbN$:
$$xE_0y\iff\exists n\,\forall m{\geq}n\,\,(x(m)=y(m)).$$
This turns out to be the smallest, in the sense of $\leq_B$,
non-smooth countable Borel equivalence relation.  This is a particular
instance of the general Glimm-Effros Dichotomy proved in
Harrington-Kechris-Louveau \cite{HKL}, but this special case can
already be derived from Effros \cite{E}.
\medskip

(iv) (Glimm-Effros Dichotomy) If $E$ is a countable Borel equivalence
relation which is not smooth, then $E_0\leq_BE$.
\medskip

(v) $E_0<_BE_\infty$.

Thus we have
$$(\identrel\bbR) <_BE_0<_BE_\infty$$
and every other countable Borel equivalence relation on an
uncountable space is in the
interval $(E_0, E_\infty )$.
\medskip

(vi) (Adams-Kechris \cite{AK}) There are continuum many pairwise
incomparable, under $\leq_B$, countable Borel equivalence relations.
\medskip

We now have all the ingredients to formulate a precise conjecture, in
response to the question about the Borel cardinality of $\mathcal D$.
This was originally formulated (as a question) in Kechris \cite {K2} and
listed (as a conjecture) in Slaman's list of Questions in Recursion
Theory, item 2.3, posted in http://math.berkeley.edu/$\sim$slaman/.
\medskip

\noindent{\bf Conjecture}: $\equiv_T$ is a universal countable Borel
equivalence relation, i.e., $(\equiv_T)\sim_BE_\infty$.

\section{Known results and implications}
\label{sec:known}

There is some information already available about the complexity of
$\equiv_T$.

\begin{theorem}
(Slaman-Steel \cite{SS}) $E_0<_B(\equiv_T)$.\end{theorem}

This has been strengthened in Kechris \cite{K1} to show that $\equiv_T$
is not amenable and in Jackson-Kechris-Louveau \cite{JKL} to show that
$\equiv_T$ is not treeable, all indications that $\equiv_T$ is quite
complex.

One of the intriguing implications of the conjecture that $\equiv_T$ is
universal concerns the existence of unusual functions on the Turing degrees.
Recall that we call a function $f:{\mathcal D}^n\rightarrow{\mathcal D}$
{\bf Borel} if there is a Borel function $F: (2^\bbN )^n\rightarrow 2^\bbN$
such that
$$f([x_1]_T,
\dots , [x_n]_T) = [F(x_1,\dots ,x_n)]_T $$
for all $x_1,\dots,x_n \in 2^\bbN$,
where $[x]_T$ is the Turing degree of $x\in 2^\bbN$.  A {\bf pairing
function}
on $\mathcal D$ is a bijection $\langle,\rangle:
{\mathcal D}^2\rightarrow{\mathcal D}$.
\medskip

\noindent{\bf Fact}.  If $\equiv_T$ is universal, then there is a Borel
pairing function on $\mathcal D$.

\begin{proof}
If $E,F$ are Borel equivalence relations on $X,Y$ respectively,
let $E\times F$ be the Borel equivalence relation on $X \times Y$
given by
$$(x,y)(E\times F)(x',y') \iff xEx' \ \&\ yFy'.$$
Clearly $E_\infty\times E_\infty\ge_BE_\infty$, so,
since $E_\infty$ is universal,
$E_\infty\times E_\infty\sim_BE_\infty$.  Hence,
if $(\equiv_T)\sim_BE_\infty$, we have
$$(\equiv_T)\times(\equiv_T)\sim_B (\equiv_T),$$
which shows that there is a Borel pairing function on
$\mathcal D$.\end{proof}

The well-known Martin Conjecture (or the 5th Victoria Delfino problem),
see Kechris-Moschovakis, Eds. \cite{KM} or Slaman's list, item 2.2,
seeks to classify definable functions on $\mathcal D$, asymptotically, i.e.,
up to identification on a cone of degrees.  One part of the conjecture
asserts, in particular, that if a Borel $f:{\mathcal D}\rightarrow{\mathcal
D}$
is not constant on a cone, then $f(d)\geq d$ on a cone.  We can now
easily see the following:
\medskip

\noindent{\bf Fact}.  If $\equiv_T$ is universal, then Martin's Conjecture
fails.

\begin{proof}
Fix $d_0\neq d_1$ in ${\mathcal D}$ and let $\langle,\rangle$
be a Borel pairing function on $\mathcal D$.  Let $f_0(d)=\langle d_0, d
\rangle$ and $f_1(d)=\langle d_1,d\rangle$.  Then $f_i:{\mathcal
D}\rightarrow
{\mathcal D}$ is Borel for $i=0,1$ and, if $A_i=\operatorname{rng}
(f_i)$, then $A_0\cap A_1=\nullset$.  Since $\equiv_T$ is countable,
one can show that the inverse of the pairing function $\langle,\rangle$
is also Borel, so the sets $A_i$ are Borel.

Clearly $f_0$ and $f_1$ are injective, so they are not constant on a cone.
Thus, if Martin's Conjecture were true, we would have that $f_i(d)\geq d$ on
a
cone for $i=0,1$.  Then $A_0$ and $A_1$ would be cofinal in the Turing
degrees,
so, by Borel Determinacy, each would contain a cone, contradiction.
\end{proof}

\section{Some more questions and answers}
\label{sec:more}

There are of course several other notions of equivalence and degree
studied in recursion theory, and similar questions and conjecture
can be considered for them too.  We will concentrate here on one of
the finest, {\it recursive isomorphism}, and one of the coarsest,
{\it arithmetic equivalence}.

Let $S_\infty$ be the group of permutations of $\bbN$,
and let $S_r$ be the subgroup consisting of all
recursive permutations.  We let
$\equiv_r$
denote {\bf recursive isomorphism} for subsets of $\bbN$.
Via our identification
of ${\mathcal P}(\bbN)$ with $2^\bbN$, we have for
$x,y\in 2^\bbN$:
$$x\equiv_ry\iff\exists \pi{\in}S_r\,\,(x\circ\pi =y).$$

For any $n\in\{ 2,3,4,\dots\}\cup\{\bbN\}$ we also define recursive
isomorphism on $n^\bbN$ by
$$x\equiv^n_ry\iff\exists\pi{\in}S_r\,\,(x\circ\pi =y),$$
so that $(\equiv^2_r)=(\equiv_r)$.

It is well-known that $(\equiv_T)\le_B(\equiv_r)$, because
$x \equiv_T y \iff x' \equiv_r y'$, where $x'$ is the
Turing jump of~$x$.  Hence, if $\equiv_T$ is universal, then
$\equiv_r$ is universal; and proving that $\equiv_r$ is universal
could be viewed as providing additional evidence that
$\equiv_T$ is universal.

Finally, we denote by $\equiv_A$ the notion of {\bf arithmetic equivalence}
on $2^\bbN$.  So $(\equiv_r)\subseteq (\equiv_T)\subseteq (\equiv_A)$.

Again, one can conjecture that $\equiv_r$ and $\equiv_A$ are universal.
Here, though, we have some answers.

\begin{theorem}
(Slaman-Steel, unpublished).  Arithmetic equivalence, $\equiv_A$,
is universal, i.e., $(\equiv_A)\sim_BE_\infty$.
\end{theorem}

So arithmetical equivalence has a Borel pairing function, and
the arithmetical analogue of Martin's Conjecture fails.

The problem for recursive equivalence is still open, but there has been a
lot
of progress.

\begin{theorem}
\label{thm:dk}
(Dougherty-Kechris \cite{DK}).  Recursive isomorphism on $\bbN^\bbN$
is universal, i.e., $(\equiv^\bbN_r)\sim_BE_\infty$.
\end{theorem}

This was very recently improved to

\begin{theorem}
\label{thm:ach}
(Andretta-Camerlo-Hjorth \cite{ACH}).  Recursive isomorphism on $5^\bbN$
is universal, i.e., $(\equiv^5_r)\sim_BE_\infty$.
\end{theorem}

However, it is not yet clear how to reduce 5 to 2.

Actually, Theorems \ref{thm:dk} and \ref{thm:ach}
are much more general.  In each case, one
actually shows that there is a fixed subgroup $S_0$ consisting of
primitive recursive (in fact much simpler) permutations such that the result
is true if $S_r$ is replaced by any countable group $S$ with $S_0\subseteq
S\subseteq S_\infty$.

There is one last problem related to Theorem~\ref{thm:dk},
that has further interesting implications.

First recall that an action of a group $G$ on a set $X$ is called {\bf free}
if $g\cdot x\neq x$ for any $x\in X$ and $g\neq 1_G$.  Also recall from
\S\ref{sec:known}
that every countable Borel equivalence relation is induced by a Borel
action of a countable group $G$.  From considerations in ergodic theory,
it turns out that it is not always possible to find a {\bf free} such
action that induces it; see Adams \cite{A}.  It has been observed though
that every known example of a countable Borel equivalence relation $E$,
which cannot be induced by a free Borel action of a countable group,
admits an invariant Borel probability measure ({\bf measure} for short).
(A measure is {\bf invariant} for $E$ if it is invariant for any Borel
action of a countable group that generates it.)  It has in fact been
conjectured that this is always the case.  In other words, a
countable Borel equivalence relation which does not admit an invariant
measure can be induced by a free Borel action of a countable group.

By using the arguments in \S 2 of Dougherty-Jackson-Kechris \cite{DJK}
and a theorem of Nadkarni \cite{N}, it can be seen that this last
assertion is equivalent to the following:
\medskip

(\dag) There is a universal countable Borel equivalence relation, which is
induced by a free Borel action of a countable group.
\medskip

We return now to Theorem~\ref{thm:dk}.
We have that $\equiv^\bbN_r$ is induced by the
following Borel action of $S_r$ on $\bbN^\bbN$:
$$\pi\cdot x=x\circ\pi^{-1}.$$
This action is not free, but its restriction to
$$[\bbN]^\bbN=\{x\in\bbN^\bbN: x\text{ is one-to-one}\}$$
is.  It is natural to conjecture that Theorem~\ref{thm:dk}
can be strengthened to the
statement that $(\equiv_r)\restrict [\bbN]^\bbN$ is universal.  If this
turns
out to be the case, this will also prove (\dag).

\section{Some proofs}

We will give here our proof of Theorem~\ref{thm:dk}
(and a related result).  This comes
from the unpublished Dougherty-Kechris \cite{DK}.  Although
Theorem~\ref{thm:dk} has now been superseded by Theorem~\ref{thm:ach},
our proof uses different methods and may find other
applications in the future.

As we indicated in \S\ref{sec:more},
one has in fact a stronger result.  For any
subgroup~$S$ of~$S_\infty$, and any $X$, let for $x,y\in X^\bbN$:
$$x\equiv^X_Sy\iff\exists\pi{\in}S\,\,(x\circ\pi =y).$$
So $(\equiv^\bbN_r)=(\equiv^\bbN_{S_r})$.  We call $S$ {\bf primitive
recursive} if $S=\{g_n: n \in\bbN\}$, with $g(n,m)=g_n(m)$ primitive
recursive.
We now have:

\begin{theorem}
\label{thm:main1}
There is a primitive recursive countable group $S_0\subseteq S_\infty$
such that for any countable group $S$ with $S_0\subseteq S\subseteq
S_\infty$,
we have that $\equiv^\bbN_S$ is a universal countable Borel equivalence
relation.  In particular this is true for $\equiv^\bbN_r$.
\end{theorem}

\begin{proof}
To explain the basic idea, consider a countable infinite
group $H$ and fix a one-to-one enumeration $H=\{h_n:n\in\bbN\}$ of it.
Then any
$h_a\in H$ corresponds to a permutation $\tilde a\in S_\infty$ given
by $h_{\tilde a(n)}=h_nh_a$ (the right regular representation).  Fix
also a bijection $\langle,\rangle :\bbN^2\rightarrow\bbN$ and let
$\pi_a\in S_\infty$ be defined by
$$\pi_a (\langle n,m\rangle )=\langle \tilde a(n), m\rangle .$$

Now given an action $(h,x)\mapsto h\cdot x$ of $H$ into a space of the
form $X^\bbN$ and the corresponding equivalence relation $E_H$, define
the function $f:X^\bbN\rightarrow X^\bbN$ by
$$f(x)(\langle n,m\rangle) =(h_n\cdot x)(m).$$
Then we have
\begin{align*}
f(h_a\cdot x)(\langle n,m\rangle )&=(h_n\cdot(h_a\cdot x))(m)\\
&=(h_{\tilde a(n)}\cdot x)(m)\\
&=f(x)(\langle\tilde a(n),m\rangle )\\
&=(f(x)\circ\pi_a)(\langle n,m\rangle );
\end{align*}
hence, $f(h_a\cdot x)=f(x)\circ \pi_a$.  It follows that if $H_0=\{\pi_a:
a\in\bbN\}$ (a countable subgroup of $S_\infty$), then
\begin{gather}
xE_Hy\iff f(x)\equiv^X_{H_0}f(y).\tag{*}
\end{gather}
Unfortunately, if $S_\infty\supseteq H'\supseteq H_0\ , H'$ a countable
group, then we cannot, in general, replace $H_0$ by $H'$ in (*) since
it could be that $f(x)\equiv^X_{H'} f(y)$ via some $\pi\in H'\setminus
H_0$.  After appropriately choosing $H$, $X$, and the action of~$H$ on
$X^\bbN$ (so that at least $E_H$ is universal), we will modify $f(x)$
to $f^*(x)\in (X^*)^\bbN$, for some $X^*$, by encoding in it some
further information, so that even if $f(x)\equiv^{X^*}_{H'} f(y)$ via some
$\pi\in H'\setminus H_0$ we can still conclude that $xE_Hy$.  In particular,
although the $X$ we will start with will be finite, this encoding will
require $X^*$ to be infinite.  Moreover, we will be forced to restrict the
$x$'s to some subset of $X^\bbN$, say $Y\subseteq X^\bbN$, so we will
also need to make sure that $E_H\restrict Y$ is universal.

We will now implement this idea.  We fix some notation first:

For any $X$ and countable group $G$, we have the shift action of $G$
on $X^G$ given by
$$g\cdot x(h)=x(g^{-1}h).$$
This induces for any subgroup $H\subseteq G$ an action of $H$ on $X^G$
and we denote the corresponding equivalence relation by $E(H,X^G)$.
If $G$ is infinite, fixing a one-to-one enumeration of $G$, we can view this
as an action of $H$ on $X^\bbN$.

Now fix a one-to-one enumeration $\{g_n:n \in \bbN\}$
of the free group $F_2$ on two generators, with $g_0=1$
where $1$~is the identity element of~$F_2$.  Define $\tilde a$ and
$\pi_a$ as above by the formulas $g_{\tilde a(n)} = g_ng_a$ and
$\pi_a(\langle n,m\rangle) = \langle \tilde a(n),m\rangle$, and let
$$S_0=\{\pi_a:a\in\bbN\}.$$
If $\{g_n: n \in \bbN\}$ and $\langle,\rangle$ are chosen appropriately,
then $S_0$ is primitive
recursive.  Fix also any countable group~$S$ such that
$S_\infty\supseteq S\supseteq S_0$;
we will show that $\equiv^\bbN_S$ is universal.  Say
$S=\{\rho_i:i\in\bbN\}$.

We call $i\in\bbN$ {\bf bad} if

(i) $\forall n\forall m\exists n'\,(\rho_i(\langle n,m\rangle )=\langle n',
m\rangle )$; and

(ii) if $\rho_i(\langle 0,m\rangle )=
\langle n_m,m\rangle$ for all $m$,
then $n_m\rightarrow \infty$ as $m \rightarrow\infty$.

\noindent
We can now easily define $n^{(i)}_j,m^{(i)}_j\in\bbN$ for $i,j\in\bbN$
such that:

(a) $0<n^{(i)}_j<n^{(i)}_{j+1}$ and $0<m^{(i)}_j<m^{(i)}_{j+1}$;

(b) $(i,j)\neq (i',j')\implies m^{(i)}_j\neq m^{(i')}_{j'}$;

(c) if $i$ is bad, then $n_{m^{(i)}_j}=n^{(i)}_j$.

Also, for the free group $F_k$ with $k$ generators and $g\in F_k$,
$m\in\bbN$,
let $B_k(g,m)$ be the ball of radius $m$ around $g$ in the tree of $F_k$;
i.e., $B_k(g,m)$ is the set of all products $gh$ where $h$ is a
word in~$F_k$ of length at most~$m$.

Now consider the shift action of $F_2$ on $9^{F_3}$ (9 is a large
enough number here) and the Borel set $A\subseteq 9^{F_3}$ defined by
\begin{align*}
y\in A\iff \forall i\forall j \,&\bigl[[(g_{n_j^{(i)}}\cdot y)\restrict
B_3(1,m^{(i)}_j)=
(g_{n^{(i)}_{j+1}}\cdot y)\restrict B_3(1,m^{(i)}_j)]\implies\\
&\quad g_{n^{(i)}_j}\cdot y=g_{n^{(i)}_{j+1}}\cdot y\bigr],
\end{align*}
where $1$~is the identity element of~$F_3$.

\begin{lemma}
$E(F_2,9^{F_3})\restrict A\leq_B(\equiv^\bbN_S)$.
\end{lemma}

\begin{proof}
Fix an injection~$c$ from the countable set $\bigcup_m 9^{B_3(1,m)}$
to~$\bbN$.  Now define $f^*: A \rightarrow \bbN^\bbN$ by
$f^*(x)=x^*$, where $x^*(\langle n,m
\rangle )=c((g_n\cdot x)\restrict B_3(1,m))$.  Thus $x^*(\langle n,m
\rangle )$ encodes the values of $g_n\cdot x$ at the ball of radius $m$
around $1\in F_3$.  In particular, $x^*(\langle n,m\rangle )$ encodes
(i.e., uniquely determines) $m$ as well.
(If we were to take $f(x)$ as in the intuitive explanation in the
beginning of this proof, then $f(x)(\langle n,m\rangle )$ would be just
$g_n\cdot x(p_m)$, where $\{p_m:m\in\bbN\}$ is a one-to-one enumeration
of~$F_3$.)

We claim that
$$xE(F_2,9^{F_3})y\iff x^*\equiv^\bbN_Sy^*,$$
which completes the proof.

$\Rightarrow$:  Clearly $y=g_a\cdot x\implies y^*=x^*\circ\pi_a$.

$\Leftarrow$: Say now $\pi\in S$ is such that $y^*=x^*\circ\pi$, i.e.,
$y^* (\langle n,m\rangle )=x^*(\pi (\langle n,m\rangle ))$.  Since
$x^* (\langle n,m\rangle )$ encodes $m$, it follows that there is a
function $\pi':\bbN\rightarrow\bbN$ such that $\pi (\langle
n,m\rangle)=\langle \pi'(\langle n,m\rangle), m \rangle$ for all
$n$ and~$m$; that is, the second
coordinate is left fixed by~$\pi$.  (Note that all $\pi_a$ have this
property, of course.  By our encoding we have forced any $\pi$ as above
to have it as well.)

We now have two cases:

(I)  $\pi'(0,m)$ does not tend to~$\infty$ as $m \rightarrow \infty$.
So there must exist a number~$\ell$ such that, for infinitely many~$m$,
$\pi'(0,m)=\ell$.  For any such~$m$, we have $y^*(\langle 0,m\rangle )
=x^*(\langle\ell, m\rangle )$, i.e., $y\restrict B_3(1,m)=(g_\ell \cdot
x)\restrict B_3(1,m)$; since there are arbitrarily large such~$m$,
it follows that $y=g_\ell \cdot x$, so
$xE(F_2,9^{F_3})y$.

(II) $\pi'(0,m)\rightarrow\infty$ as $m \rightarrow\infty$.
So if $\pi =\rho_i$, then $i$ is bad.  For any~$j$, we have
$y^*(\langle 0,m^{(i)}_j\rangle )=x^*(\langle n^{(i)}_j,m^{(i)}_j\rangle )$,
i.e., $y\restrict B_3(1,m^{(i)}_j)=(g_{n^{(i)}_j}\cdot x)\restrict
B_3(1,m^{(i)}_j)$; but we also have
$y\restrict B_3(1,m^{(i)}_{j+1})=(g_{n^{(i)}_{j+1}}\cdot x)\restrict
B_3(1,m^{(i)}_{j+1})$, and $m^{(i)}_j<m^{(i)}_{j+1}$, so we get
$(g_{n^{(i)}_j}\cdot x)\restrict B_3(1,m^{(i)}_j)=(g_{n^{(i)}_{j+1}}\cdot
x)\restrict
B_3(1,m^{(i)}_j)$.  So,
since $x\in A$, we have $g_{n^{(i)}_j}\cdot x=g_{n^{(i)}_{j+1}}\cdot x$
for all~$j$,
i.e., $g_{n^{(i)}_0}\cdot
x=g_{n^{(i)}_1}\cdot x=g_{n^{(i)}_2}\cdot x=\dotsb$.  It follows that
$y\restrict B_3(1,m^{(i)}_j)=(g_{n^{(i)}_0}\cdot x)\restrict
B_3(1,m^{(i)}_j)$ for all~$j$; since $m^{(i)}_j\rightarrow\infty$ as
$j \rightarrow \infty$, we have $y=g_{n^{(i)}_0}\cdot
x$, so $xE(F_2,9^{F_3})y$ again.
\end{proof}

It remains to show that $E(F_2,9^{F_3})\restrict A$ is universal.  For that
we will show that
$$E(F_2,2^{F_2})\leq_BE(F_2,9^{F_3})\restrict A,$$
which is enough, since $E(F_2,2^{F_2})$ is universal (see, e.g.,
Dougherty-Jackson-Kechris \cite{DJK}).

\begin{lemma}
\label{lem:2to9}
There is a Borel injection $f: 2^{F_2}\rightarrow 9^{F_3}$ with $f(2^{F_2})
\subseteq A$ which preserves the group action of~$F_2$ (i.e., for all
$g\in F_2$ and $x \in 2^{F_2}$, $f(g\cdot x)=g\cdot f(x)$).  So in
particular
$$E(F_2,2^{F_2})\leq E(F_2,9^{F_3})\restrict A.$$
\end{lemma}

To prove this lemma, we will need the following technical sublemma.

\begin{sublemma}
For each $w\in F_2\setminus \{1\}$, there is a Borel injection
$f_w: 2^{F_2}\rightarrow 6^{F_2}$ which preserves the group action
of~$F_2$ and satisfies
$$f_w(x)(g)=f_w(x)(gw)\implies g^{-1}\cdot x
=w^{-1}g^{-1}\cdot x$$
for all $g \in F_2$ and $x \in 2^{F_2}$.
\end{sublemma}

We will assume this and complete the proof.

\begin{proof}[Proof of Lemma~\ref{lem:2to9}]
Let $\{\alpha_1,\alpha_2\}$ be the generators of $F_2$ and $\{\alpha_1,
\alpha_2,\alpha_3\}$ the generators of $F_3$.
Define $f(x)$ for $x\in 2^{F_2}$ as
follows:

(i) If $g\in F_2$, then $f(x)(g)=x(g)$.

(ii) If $g=h\alpha^{-p}_3g'$, with $h\in F_2$, $p>0$, and $g'$ not starting
with
$\alpha^{\pm 1}_3$, then $f(x)(g)=2$.

(iii) If $g=h\alpha^p_3g'$, with $h,g'$ as in (ii) and $p>0$, $p\neq
m^{(i)}_j$
for all $i,j$,
then $f(x)(g)=2$.

(iv) If $g=h\alpha^{m^{(i)}_j}_3g'$, with $h,g'$ as in (ii), then $f(x)(g)=
f_{w^{(i)}_j}(x)(h)+3$, where
$$w^{(i)}_j=g_{n^{(i)}_j}g_{n^{(i)}_{j+1}}^{-1}.$$

It is easy to check that $f$ is one-to-one and preserves the action
of~$F_2$.
So it remains
to verify that $f(x)\in A$.

So fix $i,j$ with
$$(g_{n^{(i)}_j}\cdot f(x))\restrict
B_3(1,m^{(i)}_j)=(g_{n^{(i)}_{j+1}}\cdot
f(x))\restrict B_3(1,m^{(i)}_j).$$
If $d=\alpha^{m^{(i)}_j}_3$, then $d\in B_3(1,m^{(i)}_j)$, so
$$f(x)(g^{-1}_{n^{(i)}_j}d)=f(x)(g^{-1}_{n^{(i)}_{j+1}}d),$$
thus
\begin{align*}
f_{w^{(i)}_j}(x)(g^{-1}_{n^{(i)}_j})&
 =f_{w^{(i)}_j}(x)(g^{-1}_{n^{(i)}_{j+1}})\\
&=f_{w^{(i)}_j}(x)(g^{-1}_{n^{(i)}_j}w^{(i)}_j).
\end{align*}
By the sublemma, $g_{n^{(i)}_j}\cdot x
=(w^{(i)}_j)^{-1}g_{n^{(i)}_j}\cdot x = g_{n^{(i)}_{j+1}}\cdot x$, so
\begin{align*}
g_{n^{(i)}_j}\cdot f(x) &=f(g_{n^{(i)}_j}\cdot x)\\
&=f(g_{n^{(i)}_{j+1}}\cdot x)\\
&=g_{n^{(i)}_{j+1}}\cdot f(x);
\end{align*}
since $i,j$ were arbitrary, $f(x)\in A$.
\end{proof}

It remains to prove the sublemma.

\begin{proof}[Proof of Sublemma]
View $F_2$ as a rooted tree in the usual way (1 is
the root of this tree, and there is an edge between $g$ and~$g\alpha_i$
for any group element~$g$ and generator~$\alpha_i$).
Thus $x\in 2^{F_2}$ is a labeling of this tree
using labels 0,1.  Similarly for $6^{F_2}$.  Then $g^{-1}\cdot x$ is the
same
labeling except that the root of the tree is at~$g$ instead of~$1$.  So
the condition
$$\forall g\,\,[g^{-1}\cdot x\neq w^{-1}g^{-1}\cdot x\implies f_w(x)(g)
\neq f_w(x)(gw)]$$
just means that if $x$, viewed from root $g$, is different from $x$ viewed
from $gw$, then the label of $f_w(x)$ at $g$ is different from the label
of $f_w(x)$ at $gw$.  Moreover, to guarantee that $f_w(g'\cdot x)=g'
\cdot f_w(x)$ for each $g'\in F_2$, we will make sure that the value of
$f_w(x)$ at any $g$ depends only on the labeling $x$ viewed from root
$g$ (and not on $g$ itself).

Given $x\in 2^{F_2}$ and $g\in F_2$, we have two cases:

(I) $g^{-1}\cdot x=w^{-1}g^{-1}\cdot x$, i.e., $x$ looks the same from root
$g$ and root $gw$ (note that this only depends on how $x$ looks from root
$g$).

Then put $f_w(x)(g)=\langle x(g),0\rangle$, where $\langle,\rangle$ is a
bijection of $2\times 3$ with $6$.

(II) $g^{-1}\cdot x\neq w^{-1}g^{-1}\cdot x$.  So $x$ looks different from
roots $g, gw$.  In particular there is a least $n=n_g(x)$ so that for
some $i,j\in\bbZ$ and $h\in F_2$ of length~$n$ we have $x(gw^ih)
\neq x(gw^jh)$.  Clearly $n_{gw^i}(x)=n_g(x)$ for any integer~$i$ (note
that $(gw^i)^{-1}\cdot x\neq w^{-1}(gw^i)^{-1}\cdot x$ as well).

The functions $p_j:B_2(1,n_g(x))\rightarrow 2$ given by
$$p_j(h)=x(gw^jh)$$
are thus not all equal.  So fix $p\in 2^{B_2(1,n_g(x))}$ with $Z=\{j\in
\bbZ :p_j=p)\neq\nullset$ and $p$ least such (in some ordering of
$2^{B_2(1,n_g(x))}$ fixed in advance).
The value of~$p$ would be the same if we started with $gw^i$
instead of~$g$; the set~$\tilde Z$ we would get from~$gw^i$ is
a translate of~$Z$ ($j \in \tilde Z$ iff $j+i \in Z$).

Also $\{j\in\bbZ :p_j\neq p\}\neq\nullset$.  If $Z$ has a largest element
$i_0$, let $f_w(x)(g)=\langle x(g),0\rangle$, if $i_0$ is even, and $f_w
(x)(g)=\langle x(g),1\rangle$, if $i_0$ is odd.  If $Z$ has no largest
element but has a least
element~$i_0$, define $f_w(x)(g)$ the same way.
Proceed similarly if $\bbZ\setminus Z$ has a least or largest element.
So assume both $Z$ and $\bbZ\setminus Z$ are unbounded in both directions.
Put
$$Z'=\{j\in Z:j+1\not\in Z\}.$$
Let finally $f_w(x)(g)=\langle x(g),0\rangle$ if $0\in Z'$, $f_w(x)(g)=
\langle x(g),1\rangle$ if $0\not\in Z'$, but the least positive element
of $Z'$ if odd, and $f_w(x)(g)=\langle x(g),2\rangle$ if this
least positive element is even.

This completes the definition of~$f$; it is straightforward to verify that
it has the desired properties.
\end{proof}

This completes the proof of Theorem~\ref{thm:main1}.
\end{proof}

We conclude with another application of these ideas.

For a countable group $G$ consider the shift action of $G$ on $X^G$.
We call $x\in X^G$ a {\bf left-free} point if for all distinct $g, g'\in G$
there exists
$h\in G$ such that $x(hg)\neq x(hg')$.  We call $x\in X^G$ a {\bf
right-free}
or just {\bf free} point, if for all distinct $g, g'\in G$ there exists
$h\in G$ such that $x(gh)
\neq x(g'h)$; equivalently, $g\cdot x\neq g'\cdot x$ for $g \ne g'$,
or simply $g\cdot x\neq x$ for all $g\neq 1_G$.  Denote by $LF$ the set
of left-free points and $F$ the set of free points.  Note that $LF$ and~$F$
are Borel $G$-invariant subsets of $X^G$.  If $G$ is abelian, clearly
$LF=F$.
But $LF$ and~$F$ are very different for free groups in the following sense.

\begin{theorem}
The equivalence relation $E(F_3,4^{F_3})\restrict LF$ is universal for
countable
Borel equivalence relations but $E(F_3, 4^{F_3})\restrict F$ is not.
\end{theorem}

\begin{proof}
The equivalence relation $E(F_3, 4^{F_3})\restrict F$ is not universal
because it is treeable; see Kechris~\cite{K2}.
For the first assertion we will show that
$E(F_2,2^{F_2})\leq_BE(F_3,4^{F_3})
\restrict LF$.

Fix a left-free point $z_0$ in $\{2,3\}^{F_2}$.
Define then $f: 2^{F_2}\rightarrow
4^{F_3}$ by:

(i) If $h\in F_2, f(x)(h)=x(h)$.

(ii) If $h\not\in F_2$, express the reduced word for~$h$ in the form
$h=h_1\alpha^{\pm 1}_3h'$ with $h'\in F_2$,
and put $f(x)(h)=z_0(h')$.

It is easy to check that $xE(F_2,2^{F_2})y\iff f(x)E(F_3,4^{F_3}) f(y)$.
It remains to verify that $f(x)\in LF$.
Let $g$ and~$g'$ be distinct elements of~$F_3$; we must find
$h \in F_3$ such that $f(x)(hg) \ne f(x)(hg')$.

Consider two cases:

(1) $g^{-1}g'\in F_2$.  Then let $p\in F_2$ be such that $z_0(p)\neq z_0
(pg^{-1}g')$, and let $h$ be such that $hg=\alpha_3p$.  Then $f(x)(hg)=f(x)
(\alpha_3p)=z_0(p)\neq z_0(pg^{-1}g')=f(x)(hgg^{-1}g')=f(x)(hg')$.

(2) $g^{-1}g'\not\in F_2$. Let $h=g^{-1}$.  Then
$$f(x)(hg)=f(x)(1)=x(1)\in\{0,1\}$$
but
$$f(x)(hg')=f(x)(g^{-1}g')=z_0(h')\in\{2,3\}$$
for some $h'\in F_2$.
\end{proof}

\bibliographystyle{amsalpha}

\end{document}